\documentclass[11pt]{amsart}    

\usepackage[pagebackref=true,colorlinks,linkcolor=blue,citecolor=magenta]{hyperref}
\usepackage{amsmath,amssymb,latexsym} 
\usepackage{graphicx}\usepackage{centernot}

\usepackage{fullpage}
\usepackage[square, numbers]{natbib}   
\usepackage{comment}

\usepackage{enumitem}

\newtheorem{theorem}{Theorem} 

\newtheorem{alphtheorem}{Theorem}

\newtheorem{corollary}{Corollary}

\newtheorem{lemma}{Lemma}
\theoremstyle{definition}

\theoremstyle{remark}

\title{Hedetniemi's Conjecture via Altermatic Number}
\author{Meysam Alishahi}
\address{M. Alishahi, 
School of Mathematical Sciences,
Shahrood University of Technology, Shahrood, Iran}
\email{meysam\_alishahi@shahroodut.ac.ir}

\author{Hossein Hajiabolhassan}
\address{H. Hajiabolhassan, 
Department of Applied Mathematics and Computer Science, Technical University of Denmark, 
DK-{\rm 2800} Lyngby, Denmark \& Department of Mathematical Sciences,
Shahid Beheshti University, P.O. Box 19839-69411, Tehran, Iran}
\email{hhaji@sbu.ac.ir}
\begin{document}
\maketitle

\begin{abstract} 
A $50$ years unsolved conjecture by 
Hedetniemi [{\it Homomorphisms of graphs and automata, \newblock {\em Thesis (Ph.D.)--University of Michigan}, 1966}]  asserts that the chromatic number of the categorical product of two graphs 
$G$ and $H$ is $\min\{\chi(G),\chi(H)\}$. The 
present authors [{\it On the chromatic number of general {K}neser hypergraphs.
\newblock {\em Journal of Combinatorial Theory, Series B}, 2015.}] introduced 
the altermatic and the strong altermatic number of graphs as two tight lower bounds for the chromatic number of graphs. 
In this work, we prove a relaxation of  Hedetniemi's conjecture in terms of strong altermatic number. Also, we present a tight lower bound for the chromatic number of the categorical product of two graphs in term of their altermatic and strong altermatic numbers.
These results enrich the family of pair graphs $\{G,H\}$ satisfying Hedetniemi's conjecture. 
\\

\end{abstract}

\section{Introduction}
A challenging and long-standing conjecture in graph theory is
Hedetniemi's conjecture~\cite{MR2615860} which asserts that
the chromatic number of the categorical product of two graphs is the minimum of
that of graphs. There are a few general
results about the chromatic number of the categorical product of graphs whose chromatic
numbers are large enough. In view of topological bounds, it has earlier been shown that Hedetniemi's conjecture holds for any
two graphs for which some of topological bound on the chromatic number is tight,
see~\cite{MR557892, MR2718285, MR2445666}. A family of graphs is tight if  Hedetniemi's conjecture holds for any two graphs of this family. 
In this paper, we enrich the topological family of tight graphs. In this regard,
we consider the {\it altermatic number} and the {\it strong altermatic number} of graphs to determine the chromatic number of
the categorical product of some graphs.
The altermatic number of general Kneser hypergraphs was introduced in~\cite{2013arXiv1302.5394A} to obtain a tight lower bound for
the chromatic number of general Kneser hypergraphs which is a substantial improvement on
the Dol'nikov-K{\v{r}}{\'{\i}}{\v{z}} lower bound~\cite{MR953021, MR1081939, MR1665335}. 

This paper is organized as follows.
In Section 2, we set up notations and terminologies. In particular,
we will be concerned with the definition of the altermatic number and the strong altermatic number of graphs and
we mention some results about them. 
In Section 3, in view of an appropriate Kneser representation
for the Mycielskian of a graph, we show that the altermatic number
behaves like the chromatic number for the Mycielskian of a graph. Precisely, we show that
the altermatic number of the Mycielskian of a
graph $G$ is at least  the altermatic number of $G$ plus one. Next,
we prove some relaxation versions of Hedetniemi's conjecture for altermatic number and strong altermatic number.
By topological methods, it has earlier been shown that Hedetniemi's conjecture holds for any
two graphs of the family of Kneser graphs, Schrijver graphs, and  the iterated Mycielskian of
any such graphs. We enrich this result to other graphs such as a large family of Kneser multigraphs,
matching graphs, and permutation graphs.
\section{Notations and Terminologies}
\subsection{Basic Preliminaries}
Hereafter, the symbol $[n]$ stands for the set $\{1,2,\ldots, n\}$.
A {\it hypergraph} $\mathcal{ H}$ is an ordered pair $(V(\mathcal{ H}),E(\mathcal{ H}))$ consisting of a nonempty set of {\it vertices} $V(\mathcal{ H})$ and a set of {\it edges} $E(\mathcal{ H})$ which is a family of nonempty subsets of $V(\mathcal{ H})$. Unless otherwise stated, we consider simple hypergraphs, i.e.,
$E(\mathcal{ H})$ is a family of distinct nonempty subsets of $V(\mathcal{ H})$.
A mapping $h:V(\mathcal{ H})\longrightarrow [k]$ is called a $k$-coloring for $\mathcal{ H}$
if the vertices of any edge $e$
of $\mathcal{ H}$ receive  at least two distinct colors, i.e., $|\{h(v):\ v\in e\}|\geq 2$.
In other words, no edge is monochromatic. The minimum $k$ such that
$\mathcal{ H}$ admits a $k$-coloring is called the
{\it chromatic number} of $\mathcal{ H}$ and denoted by $\chi(\mathcal{ H})$.
Since  for any $k$ there is no $k$-coloring for a hypergraph
with some edge of size $1$,
we define the chromatic number of such a hypergraph to be infinite.
A hypergraph $\mathcal{ H}$ is {\it $k$-uniform}, if all edges of $\mathcal{ H}$ have the same size $k$. A $2$-uniform hypergraph is called a {\it graph}.

For two graphs $G$ and $H$, a {\it homomorphism } from $G$ to $H$
is a mapping $f:V(G)\longrightarrow V(H)$ which preserves the adjacency, i.e.,
if $xy\in E(G)$, then $f(x)f(y)\in E(H)$.
For brevity, we use $G\longrightarrow H$ to denote that there is a homomorphism from
$G$ to $H$. We call $G$ and $H$ are {\it homomorphically equivalent},  denoted by $G\longleftrightarrow H$, if we have both $G\longrightarrow H$ and $H\longrightarrow G$.
It is possible to reformulate several well-known concepts in graph colorings via
graph homomorphism.
For example, the chromatic number of a graph $G$  is the minimum integer $k$ for
which there is a homomorphism from $G$ to
the complete graph $K_k$. One can see that if two graphs are homomorphically 
equivalent, then they have the same
chromatic number, circular chromatic number, and fractional chromatic number.
An {\it isomorphism} between two graphs $G$ and $H$ is a bijective 
map $f:V(G)\longrightarrow V(H)$ such that
both $f$ and $f^{-1}$ are graph homomorphism.
For brevity, we use $G\cong H$ to mention that there is an isomorphism between
$G$ and $H$. Also, if $G\cong H$, then we say $G$ and $H$ are {\it isomorphic}.

For a hypergraph $\mathcal{ H}=(V(\mathcal{ H}), E(\mathcal{ H}))$, the general Kneser graph ${\rm KG}(\mathcal{ H})$ has the set $E(H)$ as vertex set and two vertices are adjacent if the corresponding edges are disjoint.
It is a well-known result that for any graph $G$ there is a family $\mathcal{ F}_G$ of hypergraphs such that for any $\mathcal{ H}\in \mathcal{ F}_G$,
the graph ${\rm KG}(\mathcal{ H})$ is isomorphic to $G$.
Hereafter, any hypergraph $\mathcal{ H}\in \mathcal{ F}_G$ will be referred to as a {\it Kneser representation} of $G$.
One can see that a graph has various Kneser representations. 

A subset $S\subseteq [n]$ is said to be {\it $s$-stable} if
$s\leq|i-j|\leq n-s$ for any distinct $i,j\in S$.
Hereafter, for a set $A\subseteq [n]$, the hypergraphs ${A\choose k}$ and ${A\choose k}_s$ have $[n]$ as vertex set and
the edge sets consisting of all $k$-subsets and all $s$-stable $k$-subsets of $A$, respectively.
The {\it ``usual'' Kneser graph} ${\rm KG}(n,k)$ and {\it $s$-stable Kneser graph} ${\rm KG}(n,k)_{s-stab}$
have all $k$-subsets and all $s$-stable $k$-subsets of $[n]$ as vertex sets, respectively.
Also, in these graphs, two vertices are adjacent,  whenever the corresponding sets are disjoint.
Note that ${\rm KG}({[n]\choose k})\cong{\rm KG}(n,k)$ and ${\rm KG}({[n]\choose k}_s)\cong{\rm KG}(n,k)_{s-stab}$.
The graph ${\rm KG}(n,k)_{2-stab}={\rm SG}(n,k)$ is known as {\it Schrijver graph}.
In 1955, Kneser conjectured that $\chi({\rm KG}(n,k))=n-2k+2$.
Lov\'asz \cite{MR514625}, by using the Borsuk-Ulam theorem, proved this conjecture.
Next, it was improved by Schrijver~\cite{MR512648} who proved that the
Schrijver graph ${\rm SG}(n,k)$ is a critical subgraph of the Kneser graph
${\rm KG}(n,k)$ with the same chromatic number. 
For more about the chromatic number of $s$-stable Kneser graphs, we refer the reader to~\cite{2013arXiv1302.5394A,MR2448565,MR2793613,MR512648}.

\subsection{Altermatic Number}
Let $V=\{v_1,v_2,\ldots,v_n\}$ be a set of size $n$. 
The signed-power set of $V$ is defined as $P_s(V)=\left\{(A,B)\ :\ A,B\subseteq V, A\cap B=\varnothing\right\}$.
Let $L_V$ be the set of all linear orderings of $V$, i.e.,
$L_V=\left\{v_{i_1}<v_{i_2}<\cdots<v_{i_n}\ :\ (i_1,i_2,\ldots,i_n)\in S_n\right\}$, where $S_n$ is the symmetric group. 
For any linear ordering $\sigma: v_{i_1}<v_{i_2}<\cdots<v_{i_n}\in L_V$ and $1\leq j\leq n$, define
$\sigma(j)=v_{i_j}$.
Also, for any $X=(x_1,\ldots,x_n)\in\{-1,0,+1\}^{n}$,
define $X_\sigma=(X^+_\sigma,X^-_\sigma)\in P_s(V)$, where $X^+_\sigma=\{\sigma(j):\ x_j=+1\}=\{v_{i_j}:\ x_j=+1\}$ and
$X^-_\sigma=\{\sigma(k):\ x_k=-1\}=\{v_{i_k}:\ x_k=-1\}$.

The sequence $y_1,y_2,\ldots,y_m\in \{-1,+1\}$ is said to be an {\it alternating sequence},
if any two consecutive terms of this sequence are different, i.e., $y_iy_{i+1}<0$,  for $i=1,2,\ldots,m-1$.
The {\it alternation number} $alt(X)$ of an $X=(x_1,x_2,\ldots,x_n)\in\{-1,0,+1\}^n\setminus\{(0,0,\ldots,0)\}$ is the
length of a longest alternating subsequence of nonzero entries of
$(x_1,x_1,\ldots,x_n)$. Note that we consider just nonzero entries
to determine the alternation number of $X$.
We also define $alt((0,0,\ldots,0))=0$.

For any hypergraph $\mathcal{ H}=(V, E)$ and a linear ordering $\sigma\in L_V$, define $alt_\sigma(\mathcal{ H})$ (resp. $salt_\sigma(\mathcal{ H})$)
to be the largest integer $k$ such that there is an
$X\in\{-1,0,+1\}^n$ with  $alt(X)=k$
and that none (resp. at most one) of $X^+_\sigma$ and $X^-_\sigma$ contains any (resp. some) edge of $\mathcal{ H}$.

Note that if every singleton is an edge of $\mathcal{ H}$, then
$alt_\sigma(\mathcal{ H})=0$. Also, $alt_\sigma(\mathcal{ H})\leq salt_\sigma(\mathcal{ H})$ and equality can hold.
Now set $alt(\mathcal{ H})=\min\{alt_\sigma(\mathcal{ H}):\ \sigma\in L_V\}$
and $salt(\mathcal{ H})=\min\{salt_\sigma(\mathcal{ H}):\ \sigma\in L_V\}$.
Define the {\it altermatic number}
and the {\it strong altermatic number} of a graph $G$, respectively, as follows
$$\zeta(G)=\displaystyle\max_\mathcal{ H}\left\{|V(\mathcal{ H})|-alt(\mathcal{ H}):\ {\rm KG}(\mathcal{ H})\longleftrightarrow G\right\}$$
and
$$\zeta_s(G)=\displaystyle\max_\mathcal{ H}\left\{|V(\mathcal{ H})|+1-salt(\mathcal{ H}):\ {\rm KG}(\mathcal{ H})\longleftrightarrow G\right\}.$$

It was proved in \cite{2013arXiv1302.5394A} that both altermatic number
and  strong altermatic number provide tight lower bounds for the chromatic number of graphs. 
\begin{alphtheorem}{\rm \cite{2013arXiv1302.5394A}}\label{lowermain}
For any graph $G$, we have
$$\chi(G)\geq \max\left\{\zeta(G),\zeta_s(G)\right\}.$$
\end{alphtheorem}
In view of the aforementioned theorem, one can determine the chromatic number of some family of graphs. We should 
mention that Meunier~\cite{2013arXiv1306.1112M} proved that it is a hard problem to compute the altermatic number of graphs. 

\section{Results}
In this section, first we study the altermatic number of the Mycielski construction of graphs.
Next, we present some lower bounds for the chromatic number of the categorical product of graphs in terms of the altermatic number and
strong altermatic number of graphs.
These bounds enable us to determine the chromatic number of the categorical product of some new families of graphs.

\subsection{Mycielski Construction}
For a given graph G with the vertex set $V(G)=\{u_1,\ldots,u_n\}$, the {\it Mycielskian} of $G$ is the graph $M(G)$
with the vertex set $V(M(G))=\{u_1,\ldots,u_n,v_1,\ldots,v_n,w\}$ and the edge set
$E(M(G))=E(G)\cup\{u_iv_j: u_iu_j\in E(G)\}\cup\{wv_i: 1\leq i \leq n\}.$
 In fact, the {\it Mycielski} graph $M(G)$ contains the graph $G$ itself as an isomorphic induced subgraph,
 along with $n+1$ additional vertices. For any $1\leq i \leq n$, $v_i$ is called the twin of $u_i$ and they have the same neighborhood in
$G$ and also the vertex $w$ is termed the root vertex. The vertex $w$ and $v_i$'s form a star graph.
Some coloring properties of  the Mycielski graph $M(G)$ have been studied in the literature. For instance, it is known that
$\chi(M(G))=\chi(G)+1$ and $\chi_f(M(G))=\chi_f(G)+{1\over \chi_f(G)}$, where $\chi_f(G)$ is the {\it fractional chromatic number} of $G$.

For any vector $\vec{r} = (r_1, \ldots , r_n)$ with positive integer entries and any ordering $\mu:~u_{k_1}< u_{k_2}<\cdots<u_{k_n}$ of the vertex set of $G$, the $(\vec{r},\mu)$-blow up graph
$G(\vec{r},\mu)$ of $G$ is
obtained by replacing each vertex $u_{k_i}$ of $G$ by $r_i$ vertices $u_{k_i}^1,\ldots,u_{k_i}^{r_i}$, such that for
any $1\leq l \leq r_i$ and $1\leq l'\leq r_j$, $u_{k_i}^{l}$ is adjacent to $u_{k_j}^{l'}$ if
$u_{k_i}$ is adjacent to $u_{k_j}$ in $G$. In other words, the edge $u_{k_i}u_{k_j}$ is replaced
by the complete bipartite graph $K_{r_i,r_j}$. Note that
for any vector $\vec{r}$ with positive integer entries,  two graphs $G(\vec{r},\mu)$ and $G$ are homomorphically equivalent.
\begin{lemma}\label{Myc}
For any graph $G$, we have $\zeta(M(G))\geq\zeta(G)+1$.
\end{lemma}
\begin{proof}{
Consider a hypergraph  $\mathcal{ F}=([n],E(\mathcal{ F}))$ such that that ${\rm KG}(\mathcal{ F})$ is homomorphically equivalent to $G$ and that
there exists an ordering $\sigma$ of $[n]$ for which $\zeta(G)=n-alt_{\sigma}(\mathcal{ F})$. Let
$E(\mathcal{ F})=\{A_1,\ldots,A_m\}$.
Set $t=\zeta(G)$ and $\vec{r}= (r_1, \ldots , r_{2m+1})$, where
$r_1=\cdots=r_m=2t+1$, $r_{m+1}=\cdots=r_{2m}={2t+1\choose t+1}$, and $r_{2m+1}=1$. Also,
consider an ordering $\mu$ of the vertex set of $M({\rm KG}(\mathcal{ F}))$ such that the first $m$ vertices of $\mu$ form the
vertex set of ${\rm KG}(\mathcal{ F})$ and the last vertex corresponds to the root vertex of $M({\rm KG}(\mathcal{ F}))$.

In what follows, we introduce a Kneser representation for
the $(\vec{r},\mu)$-blow up of $M({\rm KG}(\mathcal{ F}))$, i.e., $M({\rm KG}(\mathcal{ F}))(\vec{r},\mu)$. Note that $M({\rm KG}(\mathcal{ F}))(\vec{r},\mu)$
and $M(G)$ are homomorphically equivalent. Define
$$V'=\{b_1, b_2, \ldots, b_{2t+1}, c_1,c_2, \ldots, c_{(2t+1)(m-1)}\},$$
$$V''=\{a_{1,1}, a_{1,2},\ldots, a_{1,(2t+1)}, a_{2,1}, a_{2,2},\ldots, a_{2,(2t+1)}, \ldots,
a_{m,1}, a_{m,2},\ldots, a_{m,2t+1}\},$$
where $V'$, $V''$, and the set $[n]$ are pairwise disjoint. Set
$V=[n]\cup V'\cup V''$ and $l={2t+1 \choose t+1}$. For any $1\leq i \leq m$ and $1\leq j \leq 2t+1$, define $A_{i,j}=A_i\cup\{a_{i,j}\}$.
Moreover, for any $1\leq i \leq m$, consider distinct sets $B_{i,1},\ldots,B_{i,l}$ such that for any $1\leq k \leq l$,
there exists a unique $(t+1)$-subset $\{b_{k_1},b_{k_2},\ldots, b_{k_{t+1}}\}$ of $\{b_1,b_2,\ldots, b_{2t+1}\}$ where
$B_{i,k}=A_i\cup \{b_{k_1},b_{k_2},\ldots, b_{k_{t+1}}\}$. Set $\mathcal{ H}=(V,E(\mathcal{ H}))$, where
$$E(\mathcal{ H})=\{A_{i,j}, B_{i,k}: 1\leq i \leq m, 1\leq j \leq 2t+1, 1\leq k \leq l\}
\cup\{V''\}.$$
One can check that ${\rm KG}(\mathcal{ H})$ provides a Kneser representation for $M({\rm KG}(\mathcal{ F}))(\vec{r},\mu)$. To see this,
one can check that $A_{ij}$'s, $B_{ij}$'s, and $V''$
are corresponding to  the copies of the vertices of ${\rm KG}(\mathcal{ F})$, their twins,
and the root vertex, respectively. Note that
$V(\mathcal{ H})=V$ and also $c_1, \ldots, c_{(m-1)(2t+1)}$
are the isolated vertices of $\mathcal{ H}$. As a benefit of using isolated vertices,
we present an ordering $\pi$ to determine the altermatic number of $M({\rm KG}(\mathcal{ F}))(\vec{r},\mu)$.
First, consider the ordering $\tau$ as follows
$$\begin{array}{c}
\quad \quad \quad a_{1,1}< c_1< a_{2,1}< c_2< \cdots < a_{m-1,1}< c_{m-1}<  a_{m,1}< b_1<\\
\quad \quad a_{1,2}< c_{m}< a_{2,2}< c_{m+1}< \cdots < a_{m-1,2}< c_{2m-2}<  a_{m,2}< b_2 <\\
\vdots\\
a_{1,2t+1} < c_{2t(m-1)+1}< a_{2,2t+1}< \cdots< c_{(2t+1)(m-1)}<  a_{m,2t+1} < b_{2t+1}
\end{array}$$
Construct the ordering $\pi$ by concatenating the ordering $\sigma$ after $\tau$, i.e., $\pi=\tau||\sigma$.
Note that the number of elements of $\pi$ is $(2t+1)m+(2t+1)+(2t+1)(m-1)+n=2m(2t+1)+n$. Define $p=2m(2t+1)+n$.
Now we claim that
$alt_{\pi}(\mathcal{ H})\leq alt_{\sigma}(\mathcal{ F})+2m(2t+1)-1$ which implies the assertion.
To see this, assume that
$X=(x_1,x_2,\ldots,x_p) \in\{-1,0,+1\}^p\setminus \{(0,0,\ldots,0)\}$ and
$alt(X)=alt_{\sigma}(\mathcal{ F})+2m(2t+1).$
We show that $X^+_\pi$ or  $X^-_\pi$ contains an edge of $\mathcal{ H}$.
If $alt((x_1,x_2,\ldots,x_{2m(2t+1)}))=2m(2t+1)$, then $X^+_\pi$ or  $X^-_\pi$
contains $\{a_{1,1}, a_{1,2},\ldots, a_{m,2t+1}\}$,
i.e., the root vertex; and consequently, the assertion follows.
Hence, let $alt((x_1,x_2,\ldots,x_{2m(2t+1)}))\leq 2m(2t+1)-1$; and
consequently, $alt(Y)\geq alt_{\sigma}(\mathcal{ F})+1$,
where $Y=(x_{2m(2t+1)+1},x_{2m(2t+1)+2},\ldots,x_p)$. Hence,
$Y^+_\sigma$ or  $Y^-_\sigma$ contains an edge of $\mathcal{ F}$.
Without loss of generality, suppose that $Y^+_\sigma$ contains $A_i\in \mathcal{ F}$. If
$a_{i,j}=+1$ for some $1\leq j\leq 2t+1$, then $A_{i,j}\in E(\mathcal{ H})$ and
$A_{i,j}\subseteq X^+_\pi$. Moreover, if
there exists a $(t+1)$-subset $\{b_{k_1},b_{k_2},\ldots, b_{k_{t+1}}\}$
such that  $b_{k_j}=+1$ for any $1\leq j\leq t+1$, then
$A_{i}\cup \{b_{k_1},b_{k_2},\ldots, b_{k_{t+1}}\}\in E(\mathcal{ H})$ and $A_{i}\cup \{b_{k_1},b_{k_2},\ldots, b_{k_{t+1}}\}\subseteq X^+_\pi$.
On the other hand, if for some $1\leq j\leq 2t+1$, $a_{i,j}\not =+1$
and $b_j\not =+1$,  then one can check that $alt((x_{2m(j-1)+1}, x_{2m(j-1)+2},\ldots,x_{2mj}))\leq 2m-1$.
Therefore, if for any $1\leq j\leq 2t+1$, $a_{i,j}\not =+1$, and
if for any $(t+1)$-subset $\{b_{k_1},b_{k_2},\ldots, b_{k_{t+1}}\}$, there exists some
$1\leq j\leq 2t+1$ such that $b_{k_j}\not =+1$, then one can conclude that $alt((x_1,x_2,\ldots,x_{2m(2t+1)})) \leq
2m(2t+1)-(t+1)=2m(2t+1)-n+alt_{\sigma}(\mathcal{ F})-1$. Also, $alt(X)=alt_{\sigma}(\mathcal{ F})+2m(2t+1)$;
hence, we should have $alt(Y)\geq n+1$  which is impossible.
}\end{proof}
A graph $G$ is said to be {\it alternatively $t$-chromatic} (resp. {\it strongly alternatively $t$-chromatic}) if
$\zeta(G)=\chi(G)=t$ (resp. $\zeta_s(G)=\chi(G)=t$).
Note that if two graphs $G$ and $H$ are homomorphically equivalent, then $M(G)$ and $M(H)$
are homomorphically equivalent as well.
In other words, the previous theorem states that the graph $M(G)$ is
alternatively $(t+1)$-chromatic graph provided that $G$ is alternatively $t$-chromatic graph.
\subsection{Chromatic Number of Categorical Product}
There are several kinds of graph products in the literature.
The {\em categorical product } $G\times H$ of two graphs $G$ and $H$ is defined by
$V(G\times H)=V(G)\times V(H)$, where two
vertices $(u_1, u_2)$ and $(v_1, v_2)$ are adjacent if $u_1v_1\in E(G)$ and
$u_2v_2\in E(H)$.
It is easy to check that by any coloring of $G$ or $H$, we can present a coloring of
$G\times H$; and therefore, $\chi(G\times H)\leq\min\{\chi(G),\chi(H)\}$.
In 1966, Hedetniemi~\cite{MR2615860} introduced his interesting conjecture,
celled  {\it Hedetniemi's conjecture},  about the chromatic
number of the categorical product of two graphs, which states that $\chi(G\times H)=\min\{\chi(G),\chi(H)\}$.
This conjecture has been studied in the literature, see \cite{Dochtermann2009490,MR1804825, MR1185007, MR2176429,
MR2445666, MR2825542}. In \cite{MR1815614}, Zhu generalized Hedetniemi's conjecture to the circular chromatic number of graphs. Precisely, he
conjectured $\chi_c(G\times H)=\min\{\chi_c(G),\chi_c(H)\}$.
For more results about this conjecture and circular chromatic number, we refer the reader to \cite{MR1905132, MR2171371, MR1815614, MR2249284}. 
In view of topological bounds for chromatic number, it was shown that Hedetniemi's conjecture holds for any
two graphs for which some of topological bounds on the chromatic number are tight,
see~\cite{Dochtermann2009490,MR557892, MR2383129, MR2718285, MR2445666}.
\begin{alphtheorem}{\rm \cite{MR2383129, MR2718285}}\label{coindB}
For any two graphs $G$ and $H$, we have
$${\rm coind}(B(G\times H))=\min \{{\rm coind}(B(G)), {\rm coind}(B(H))\}.$$
\end{alphtheorem}
Note that it is~not known whether the following inequality holds
$${\rm coind}(B_0(G\times H)) \geq \min \{{\rm coind}(B_0(G)), {\rm coind}(B_0(H))\}.$$

The quantities ${\rm coind}(B_0(G))$ and ${\rm coind}(B(G))$ used in the statement 
of the previous assertions are referring to two topological parameters;  
coindices of two variations of box complexes of $G$, for more about these parameters see~\cite{MR1988723}.   

In this section, we present some lower bounds
for the chromatic number of the categorical product of graphs.
\begin{lemma}\label{homhalt}
Let $G$ and $H$ be two graphs. If there exists a graph homomorphism $f:H\longrightarrow G$,
then $\zeta(H)\leq\zeta(G)$ and also $\zeta_s(H)\leq\zeta_s(G)$.
\end{lemma}
\begin{proof}{
We prove $\zeta(H)\leq\zeta(G)$ and similarly one can show $\zeta_s(H)\leq\zeta_s(G)$.
First, we assume that $H$ is a subgraph of $G$ and we prove a stronger assertion. In fact, we show that
if $H$ is a subgraph of $G$ and
$\mathcal{ H}=([n], E(\mathcal{ H}))$ is a Kneser representation of $H$, i.e., ${\rm KG}(\mathcal{ H})$ is isomorphic to $H$, then
there exists a Kneser representation $\mathcal{ G}= (Y,E(\mathcal{ G}))$ for $G$ such that $n-alt(\mathcal{ H})\leq  |Y|-alt(\mathcal{ G})$.

Without loss of generality, suppose that $alt(\mathcal{ H})=alt_{\sigma}(\mathcal{ H})$, where $\sigma$ is an ordering of $V(\mathcal{ H})$. Let
$g: V(H)\longrightarrow E(\mathcal{ H})$ be an isomorphism between $H$ and ${\rm KG}(\mathcal{ H})$.
Also, one can assume that $H$ is a spanning subgraph of $G$.
To see this, let $v\in V(G)\setminus V(H)$.
Now add this vertex to $H$ as an isolated vertex to obtain $H_1$.
Set
$\mathcal{ H}_1=([n+1], E(\mathcal{ H}_1))$, where $E(\mathcal{ H}_1)=E(\mathcal{ H})\cup\{\{1,2,\ldots,n+1\}\}$. One can see that
${\rm KG}(\mathcal{ H}_1)$ is isomorphic to $H_1$. Also,
$alt_{\sigma'}(\mathcal{ H}_1)\leq alt_{\sigma}(\mathcal{ H})+1$, where the ordering $\sigma'$ is obtained by adding the element $n+1$ at the end of the ordering $\sigma$.
Therefore,
$ n+1-alt(\mathcal{ H}_1)\geq n+1-alt_{\sigma'}(\mathcal{ H}_1)\geq  n-alt_{\sigma}(\mathcal{ H})=n-alt(\mathcal{ H}).$
By repeating the previous procedure, if it is necessary, we can find a spanning subgraph
$\bar{H}$ of $G$ and a Kneser representation $\bar{\mathcal{H}}$ for $\bar{H}$ such that 
$\bar{\mathcal{H}}=([n+l], E(\bar{\mathcal{H}}))$,
$l=|V(G)|-|V(H)|$, and that $n-alt(\mathcal{ H})\leq  n+l-alt(\bar{\mathcal{ H}})$.
Therefore, it is enough to prove this lemma just for spanning subgraphs.

Now assume that $H$ is a spanning subgraph of $G$.
Again, without loss of generality, suppose that there is an edge $e=ab\in E(G)$
such that $H+e=G$.
Let $g(a)=A$ and $g(b)=B$.
Since $a$ and $b$ are~not adjacent, $A\cap B$ is a nonempty set.
Let $A\cap B=\{y_1,y_2,\ldots,y_t\}$.
Consider $2t$ positive integers $\{y'_1,y'_2,\ldots,y'_t\}\cup\{z_1,z_2,\ldots,z_t\}$ disjoint from $[n]$.
Let $\mathcal{ H}_0=\mathcal{ H}$, $\sigma_0=\sigma$, $g_0=g$, and $Y_0=[n]$. For any $1\leq i\leq t$, set $Y_i=[n]\cup\{y'_1,z_1,y'_2,z_2,\ldots,y'_i,z_i\}$.
Assume that $\mathcal{ H}_i=(Y_i,E(\mathcal{ H}_i))$, $\sigma_i\in L_{Y_i}$, and $g_i: V(H)\longrightarrow \mathcal{ H}_i$ are defined when $i<t$.
Define $g_{i+1}(a)=g_{i}(a)\cup\{y'_{i+1}\}\setminus \{y_{i+1}\}$ and  $g_{i+1}(b)=g_{i}(b)$. Also, for $u\not \in \{a,b\}$,
if $y_{i+1}\in g_i(u)$, then define $g_{i+1}(u)=g_{i}(u)\cup\{y'_{i+1}\}$; otherwise, set $g_{i+1}(u)=g_i(u)$.
Consider the hypergraph  $\mathcal{ H}_{i+1}=(Y_{i+1}, E(\mathcal{ H}_{i+1}))$, where $Y_{i+1}=Y_i\cup\{y_{i+1}',z_{i+1}\}$  and $E(\mathcal{ H}_{i+1})=\{g_{i+1}(v):\ v\in V(H)\}$.
To obtain the ordering $\sigma_{i+1}$, replace $y_{i+1}$ with $y_{i+1}<z_{i+1}< y'_{i+1}$ in the ordering
$\sigma_i$, i.e., put  $z_{i+1}$ immediately after $y_{i+1}$  and put $y'_{i+1}$ after $z_{i+1}$.
Note that $\sigma_{i+1}\in L_{Y_{i+1}}$,
where $Y_{i+1}=[n]\cup\{y'_1,z_1,y'_2,z_2,\ldots,y'_{i+1},z_{i+1}\}$. 
One can see that $alt_{\sigma_{i+1}}(\mathcal{ H}_{i+1})\leq 2+alt_{\sigma_{i}}(\mathcal{ H}_{i})$;
and therefore, $alt_{\sigma_{t}}(\mathcal{ H}_{t})\leq 2t+alt_{\sigma_{0}}(\mathcal{ H}_{0})$.
Note that ${\rm KG}(\mathcal{ H}_{t})\cong H+e=G$.
Set $\mathcal{ G}=\mathcal{ H}_{t}$ and $Y=Y_t$.
Consequently,
$$|Y|-alt(\mathcal{ G})\geq n+2t-alt_{\sigma_t}(\mathcal{ G})\geq n+2t-(2t+alt_{\sigma_{0}}(\mathcal{ H}_{0}))
=n-alt_{\sigma}(\mathcal{ H})=n-alt(\mathcal{ H}).$$

Now assume that for two graphs $G$ and $H$, there exists a graph homomorphism $f:H\longrightarrow G$.
In view of definition of $\zeta(H)$ and since $\zeta(H)< \infty$ ($\zeta(H)\leq \chi(H)$), there is a graph $\bar{H}$ and a hypergraph $\bar{\mathcal{ H}}=([m],E(\bar{\mathcal{ H}}))$ such that $\bar{H}$ is homomorphically equivalent to $H$,
${\rm KG}(\bar{\mathcal{ H}})\cong \bar{H}$ and   $\zeta(H)=m-alt(\bar{\mathcal{ H}})$.
Since $H$ and $\bar{H}$ are homomorphically equivalent and that there exists a graph homomorphism $f:H\longrightarrow G$,
the graph $\bar{G}=G\cup\bar{H}$, i.e., the disjoint union of $G$ and $\bar{H}$, is homomorphically equivalent
to the graph $G$ and also this graph has $\bar{H}$ as its subgraph.
In view of the aforementioned discussion, there are $\bar{Y}$ and $\bar{\mathcal{G}}=(\bar{Y},E(\bar{\mathcal{G}}))$ such that
${\rm KG}(\bar{\mathcal{G}})\cong \bar{G}$ and
$$\zeta(H)=m-alt(\bar{\mathcal{ H}})\leq  |\bar{Y}|-alt(\bar{\mathcal{G}})\leq \zeta(\bar{G})=\zeta(G).$$
Similarly, one can show $\zeta_s(H)\leq\zeta_s(G)$.
}\end{proof}
In view of the proof of the aforementioned lemma, the next corollary follows.
\begin{corollary}
If $H$ is a subgraph of $G$, then $\zeta(H)\leq\zeta(G)$ and $\zeta_s(H)\leq\zeta_s(G)$. In particular, for
any graph $G$, we have $\min\{\zeta(G),\zeta_s(G)\}\geq \omega(G)$, where $\omega(G)$ is the clique number of
$G$.
\end{corollary}
In the next result, we show the accuracy of Hedetniemi's
conjecture for the strong altermatic number of graphs.
\begin{theorem}\label{hedet1}
For any two graphs $G$ and $H$, we have\\
{\rm a)} $\zeta_s(G\times H)=\min\{\zeta_s(G),\zeta_s(H)\},$\\
{\rm b)} $\zeta(G\times H)\geq \max\{\min\{\zeta(G),\zeta_s(H)-1\},\min\{\zeta_s(G)-1,\zeta(H)\}\}.$
\end{theorem}
\begin{proof}{
First, we prove part (a). Consider two hypergraphs $\mathcal{ G}=(V,E(\mathcal{ G}))$ and $\mathcal{ H}=(V',E(\mathcal{ H}))$ such that
${\rm KG}(\mathcal{ G})$ and ${\rm KG}(\mathcal{ H})$ are homomorphically equivalent to $G$ and $H$, respectively, and also, $\zeta_s(G)=1+|V|-salt(\mathcal{ G})$ and
$\zeta_s(H)=1+|V'|-salt(\mathcal{ H})$.
Without loss of generality, suppose that $V=\{1,2,\ldots,n\}$, $V'=\{n+1,n+2,\ldots,n+m\}$,
$salt(\mathcal{ G})=salt_{\sigma}(\mathcal{ G})$, and $salt(\mathcal{ H})=salt_{\tau}(\mathcal{ H})$, where $\sigma$ (resp. $\tau$) is an ordering of $V$
(resp. $V'$).

Let $\mathcal{ F}=([n+m], E(\mathcal{ F}))$, where $E(\mathcal{ F})=\{A\cup B:\ A\in E(\mathcal{ G})\ \&\ B\in E(\mathcal{ H}) \}$.
One can check that ${\rm KG}(\mathcal{ F})\cong {\rm KG}(\mathcal{ G})\times {\rm KG}(\mathcal{ H})\longleftrightarrow G\times H$ and
also, $\zeta_s(G\times H)\geq 1+n+m-salt(\mathcal{ F})\geq \min\{\zeta_s(G),\zeta_s(H)\}$.
To see this, it is enough to show $salt_{\pi}(\mathcal{ F})\leq \max\{|V'|+salt(\mathcal{ G}), |V|+salt(\mathcal{ H})\}$, where $\pi=\sigma||\tau$.
Define $l=\max\{|V'|+salt(\mathcal{ G}), |V|+salt(\mathcal{ H})\}$.
Consider an $X\in \{-1,0,1\}^{n+m}$ with $alt(X)\geq 1+l$.

Now consider two vectors $X(1),X(2)\in \{-1,0,1\}^{n+m}$
such that the first $n$ coordinates of $X(1)$ (resp. the last $m$
coordinates of $X(2)$) are the same as $X$ and the last $m$
coordinates of $X(1)$ (resp. the first $n$ coordinates of $X(2)$) are zero.
If $alt(X(1))> salt_{\sigma}(\mathcal{ G})$ and $alt(X(2))> salt_{\tau}(\mathcal{ H})$,
then each of $X_{\pi}^+$ and $X_{\pi}^-$ has an edge of $\mathcal{ F}$ and it completes the proof.
On the contrary, suppose $alt(X(1))\leq salt_{\sigma}(\mathcal{ G})$. Therefore,
$alt(X)\leq alt(X(1))+alt(X(2))\leq salt_{\sigma}(\mathcal{ G})+|V'|\leq l$ which is a contradiction.
By a similar argument, we conclude that $alt(X(2))> salt_{\tau}(\mathcal{ H})$. Hence,
$\zeta_s(G\times H)\geq  \min\{\zeta_s(G),\zeta_s(H)\}$.
On the other hand, there exists a graph homomorphisms from $G\times H$ into $G$ (resp. $H$). Consequently, by
Lemma~\ref{homhalt}, we have $\zeta_s(G\times H)\leq\min\{\zeta_s(G),\zeta_s(H)\}$.

Now we prove part (b). In view of  symmetry, it suffices to prove
$ \zeta(G\times H) \geq \min\displaystyle\{ \zeta(G),\zeta_s(H)-1\}.$
Consider two hypergraphs $\mathcal{ G}=(V,E(\mathcal{ G}))$ and $\mathcal{ H}=(V',E(\mathcal{ H}))$,
such that $G\longleftrightarrow{\rm KG}(\mathcal{ G})$, $H\longleftrightarrow{\rm KG}(\mathcal{ H})$,
$\zeta(G)=|V|-alt(\mathcal{ G})$, and $\zeta_s(H)=|V'|+1-salt(\mathcal{ H})$.
Without loss of generality, suppose that $V=\{1,2,\ldots,n\}$, $V'=\{n+1,n+2,\ldots,n+m\}$,
$\zeta(G)=|V|-alt_{\sigma}(\mathcal{ G})$, and $\zeta_s(H)=|V'|+1-salt_{\tau}(\mathcal{ H})$, where
 $\sigma\in L_{V}$ and $\tau\in L_{V'}$.
Let $\mathcal{ L}=([n+m], E(\mathcal{ L}))$, where
$E(\mathcal{ L})=\displaystyle\left\{A\cup B:\ A\in \mathcal{ G}\ \& \ B\in \mathcal{ H}\right\}$.
One can check that ${\rm KG}(\mathcal{ L})\cong {\rm KG}(\mathcal{ G})\times {\rm KG}(\mathcal{ H})$; and
therefore, ${\rm KG}(\mathcal{ L})\longleftrightarrow G\times H.$
Set $\pi=\sigma||\tau$ and $M=\max\{|V|+salt_\tau(\mathcal{ H}), |V'|+alt_\sigma(\mathcal{ G})\}$.
Now we show that $alt_{\pi}(\mathcal{ L})\leq M$. To see this, consider an
$X\in\{-1,0,+1\}^{m+n}\setminus\{(0,0,\ldots,0)\}$ such that $alt(X)\geq M+1$.
Let $X(1),X(2)\in \{-1,0,1\}^{m+n}$ be the same as in the proof of the previous part.
One can see that
there exists an alternative subsequence  of nonzero terms in $X(1)$ of length more
than $alt_\sigma(\mathcal{ G})$.
Therefore, $X(1)^+_{\sigma}$ or $X(1)^-_{\sigma}$
contains some edge of $\mathcal{ G}$.
Now we show that that both $X(2)^+_{\tau}$ and $X(2)^-_{\tau}$ have some edges of $\mathcal{ H}$.
On the contrary, suppose that this is~not true.
Therefore, we have $alt(X(2))\leq salt_\tau(\mathcal{ H})$; and thus, $alt(X)\leq |V|+alt(X(2))\leq M$ which is a contradiction.
Hence, assume $A\subseteq X(2)^+_{\tau}$, $B\subseteq X(2)^-_{\tau}$ and $C\subseteq X(1)^+_{\sigma}$, where $A, B\in E(\mathcal{ H})$ and $C\in E(\mathcal{ G})$.
Now in view of $A\cup C \subseteq X^+_{\pi}$ and that $A\cup C\in E(\mathcal{ L})$, the assertion follows.

Hence, we have
$$
\begin{array}{lll}
 \zeta(G\times H) & \geq  & m+n-alt_{\pi}(\mathcal{ L})\\
 			      & \geq  & m+n-M\\
                            &   =   & \min\displaystyle
                                    \{n-alt_\sigma(\mathcal{ G}),m-salt_\tau(\mathcal{ H})\}\\
                                       &   =   & \min\displaystyle
                                    \{ \zeta(G),\zeta_s(H)-1\},
\end{array}$$
as desired.
}\end{proof}
\begin{theorem}\label{hedet2}
Let $\mathcal{ G}=(V,E(\mathcal{ G}))$ and $\mathcal{ H}=(V',E(\mathcal{ H}))$ be two hypergraphs. Also, assume that
$\sigma\in L_V$ and $\tau\in L_{V'}$ such that $\zeta({\rm KG}(\mathcal{ G}))=|V|-alt_\sigma(\mathcal{ G})$
and $\zeta({\rm KG}(\mathcal{ H}))=|V'|-alt_\tau(\mathcal{ H})$.
If
$$\max\left\{|V|+alt_\tau(\mathcal{ H}), |V'|+alt_\sigma(\mathcal{ G})\right\}\geq salt_\sigma(\mathcal{ G})+salt_\tau(\mathcal{ H}),$$
 then
 $\zeta({\rm KG}(\mathcal{ G})\times {\rm KG}(\mathcal{ H}))\geq \min\left\{\zeta({\rm KG}(\mathcal{ G})),\zeta({\rm KG}(\mathcal{ H}))\right\}.$
\end{theorem}
\begin{proof}{
Let $G={\rm KG}(\mathcal{ G})$ and $H={\rm KG}(\mathcal{ H})$.
Without loss of generality, suppose that $V=\{1,2,\ldots,n\}$ and $V'=\{n+1,n+2,\ldots,n+m\}$.
Let  $\mathcal{ L}=([n+m],E(\mathcal{ L}))$ where $E(\mathcal{ L})=\displaystyle\left\{A\cup B:\ A\in \mathcal{ G}\ \& \ B\in \mathcal{ H}\right\}$.
Note that ${\rm KG}(\mathcal{ L})\cong {\rm KG}(\mathcal{ G})\times {\rm KG}(\mathcal{ H})\cong G\times H.$
Set $\pi=\sigma||\tau$ and $M=\max\{|V|+alt_\tau(\mathcal{ H}), |V'|+alt_\sigma(\mathcal{ G})\}$. In view of the assumption, we have
$$M=\max\{|V|+alt_\tau(\mathcal{ H}), |V'|+alt_\sigma(\mathcal{ G})\}\geq salt_\sigma(\mathcal{ G})+salt_\tau(\mathcal{ H}).$$
Now we show $alt_{\pi}(\mathcal{ L})\leq M$. To see this, consider an
$X\in\{-1,0,+1\}^{m+n}\setminus\{(0,0,\ldots,0)\}$ such that $alt(X)\geq M+1$.
Consider two vectors $X(1),X(2)\in \{-1,0,1\}^{m+n}$ such that
the first $n$ coordinates of $X(1)$ (resp. the last $m$
coordinates of $X(2)$) are the same as $X$ and the last $m$
coordinates of $X(1)$ (resp.
the first $n$ coordinates of $X(2)$) are zero. One can see that
there exists an alternative subsequence  of nonzero terms in $X(1)$ (resp. $X(2)$) of length more
than $alt_\sigma(\mathcal{ G})$ (resp. $alt_\tau(\mathcal{ H})$).
Therefore,  $X(1)^+_{\sigma}$ or $X(1)^-_{\sigma}$
(resp.  $X(2)^+_{\tau}$ or $X(2)^-_{\tau}$)
has some edge of $\mathcal{ G}$
(resp. $\mathcal{ H}$).
Now we show that  both $X(1)^+_{\sigma}$ and $X(1)^-_{\sigma}$  or both $X(2)^+_{\tau}$ and $X(2)^-_{\tau}$  have some edges of $\mathcal{ G}$ or $\mathcal{ H}$,
respectively.
On the contrary, suppose that this is~not true.
Therefore, we have $alt(X(1))\leq salt_\sigma(\mathcal{ G})$ and $alt(X(2))\leq salt_\tau(\mathcal{ H})$. These inequalities imply that $alt(X)\leq alt(X(1))+alt(X(2))\leq M$ which is a contradiction. Hence, without loss of generality, suppose that  $A\subseteq X(1)^+_{\sigma}$, $B\subseteq X(1)^-_{\sigma}$, and
$C\subseteq X(2)^+_{\tau}$,where
$A, B \in E(\mathcal{ G})$ and $C\in E(\mathcal{ H})$.
Now in view of $A\cup C \subseteq X^+_{\pi}$ and that $A\cup C\in E(\mathcal{ L})$, the assertion follows.

Hence, we have
$$
\begin{array}{lll}
 \zeta(G\times H) & \geq  & m+n-alt_{\pi}(\mathcal{ L})\\
 			      & \geq  & m+n-M\\
                                       &   =   & \min\displaystyle
                                    \{ \zeta(G),\zeta(H)\},
\end{array}$$ as desired.
}\end{proof}
In view of Theorems~\ref{hedet1}~and~\ref{hedet2}, one can determine the chromatic number of the categorical product of some family of graphs. In particular, if  both of them are strongly alternatively $t$-chromatic graphs. Note that the following graphs are strongly alternatively $t$-chromatic graphs.

\begin{enumerate}
\item Schrijver graphs and Kneser graphs

\item The  Kneser multigraph  ${\rm KG}(G,\mathcal{ F})$:  $G$ is a multigraph such that all of its
edges have even multiplicities and $\mathcal{ F}$ is a family of its simple subgraphs, see~\cite{2014arXiv1401.0138A}.

\item Some of the matching graphs ${\rm KG}(G,rK_2)$, see~\cite{MatchingNew}.
\end{enumerate}

By Theorem~\ref{coindB}, for any two graphs $G$ and $H$, we have $${\rm coind}(B(G\times H))=\min \{{\rm coind}(B(G)), {\rm coind}(B(H))\}.$$ Also, in~\cite{MR2455224}, it was proved that
${\rm coind}(B(M(G))) \geq {\rm coind}(B(G))+1$, where
$M(G)$ is the Mycielskian of $G$. The present authors~\cite{2016arXiv160708780A} proved that 
for a graph $G$, we have  
$$\chi(G) \geq {\rm coind}(B_0(G))+1\geq \zeta(G)\quad\mbox{and}\quad\chi(G) \geq {\rm coind}(B(G))+2\geq  \zeta_s(G).$$
Consequently, one can see that Hedetniemi's conjecture
holds for any
two graphs of the family of strongly alternatively $t$-chromatic graphs and  the iterated Mycielskian of
any such graphs.

We showed that the following graphs are alternatively $t$-chromatic graphs. Hence, in view of
Theorems~\ref{hedet1}~and~\ref{hedet2}, one can introduce several tight families.

\begin{enumerate}
\item Kneser graphs and multiple Kneser graphs: In~\cite{2013arXiv1302.5394A}, multiple Kneser graphs
were introduced as a generalization of Kneser graphs.

\item  Kneser Multigraphs, see \cite{2014arXiv1401.0138A}.

\item A family of matching graphs, see \cite{MatchingNew}.

\item The permutation graph  $S_r(m,n)$: $m$ is large enough, see \cite{MatchingNew}

\item Any number of iterations of the  Mycielski construction starting with any graph appearing on the list above.
\end{enumerate}

\noindent{\bf Acknowledgement:}
The authors would like to express their deepest gratitude to Professor Carsten~Thomassen for his insightful comments.
They also appreciate the detailed valuable comments of  Dr.~Saeed~Shaebani.
A part of this paper was written while Hossein Hajiabolhassan was visiting School of Mathematics, Institute for Research in Fundamental Sciences~(IPM). He acknowledges the support of IPM.
Moreover, they would like to thank Skype for sponsoring their endless conversations in two countries.
\def\cprime{$'$} \def\cprime{$'$}

\end{document}